\documentclass[12pt]{article}
\usepackage{amsxtra,amssymb,amsthm,amsmath,latexsym}

\theoremstyle{plain}

\def\oH{\buildrel\circ\over H}
\def\oH1{\buildrel\circ\over H\kern-.02in{}^1}

\begin{document}


\title{ Compactness of Embeddings 
   \thanks{key words: Banach spaces, compactness, embedding operator
    }
}

\author{
A.G. Ramm\\
Mathematics Department, 
Kansas State University, \\
 Manhattan, KS 66506-2602, USA\\
ramm@math.ksu.edu\\
}

\date{}

\maketitle\thispagestyle{empty}

\begin{abstract} 
An improvement of the author's result, proved in 1961, 
concerning necessary and sufficient conditions for
the compactness of embedding operators is given.
A counterexample to a  published statement concerning 
compactness of embedding operators is constructed.
\end{abstract}


\section{Introduction}

The basic result of this note is:

{\bf Theorem 1.} {\it Let $X_1\subset X_2\subset X_3$ be Banach spaces,
$||u||_1\geq ||u||_2\geq ||u||_3$ (i.e., the norms are comparable) and  
if $||u_n||_3\to 0$ as $n\to \infty$ and $u_n$ is fundamental in $X_2$, 
 then 
$||u_n||_2\to 0$,  (i.e., the norms in $X_2$ and $X_3$ are compatible).
Under the above assumptions the embedding operator $i: X_1\to X_2$ is 
compact if and only if
the following two conditions are valid:

a) The embedding operator  $j: X_1\to X_3$ is compact,

and the following inequality holds:

b) $||u||_2\leq s ||u||_1 + c(s)||u||_3, \,\,\,\forall  
u\in X_1,$ $\forall s\in (0,1)$, where $c(s)>0$ is a constant.}

This result is an improvement of the author's old result,
originally proved in 1961 (see [2]), where $X_2$ 
was assumed to be a Hilbert space. The proof of Theorem 1 is simpler than 
the one in [2]. This proof is borrowed from the recent paper [3].
In addition to this proof, we construct a counterexample to a theorem in 
[1], p.35, where
the validity of the inequality b) in Theorem 1 is claimed without 
the assumption of the compatibility of the norms of $X_2$ and $X_3$.
(see Remark 1 at the end of this note). This counterexample is new.

\section{Proof}

1. {\it The sufficiency of conditions a) and b) for compactness of 
$i: X_1\to X_2$}

Assume that a) and b) hold and let us prove the compactness of 
$i$. Let $S=\{u: u\in X_1, ||u||_1=1\}$ be the unit sphere in $X_1$.
Using assumption a), select a sequence $u_n$ which converges 
in $ X_3$. 
We claim that this sequence 
converges also in $ X_2$. Indeed, since $||u_n||_1=1$, one 
uses assumption b) to get 
$$||u_n-u_m||_2\leq
s||u_n-u_m||_1+c(s)||u_n-u_m||_3\leq 2s +c(s)||u_n-u_m||_3.$$
Let $\eta>0$ be an arbitrary small given number. Choose $s>0$ such that 
$2s<\frac 1 2\eta$, and for a fixed $s$ choose $n$ and $m$ so large that
$c(s)||u_n-u_m||_3<\frac 1 2\eta$. This is possible because the
sequence $u_n$ converges in $ X_3$. Consequently, 
 $||u_n-u_m||_2\leq \eta$ if  $n$ and $m$ are sufficiently large.
This means that the sequence $u_n$  converges in $ X_2$.
Thus, the embedding  $i: X_1\to X_2$ is compact.
In the above argument, i.e., in the proof of sufficiency, the 
compatibility of the norms was not 
used. 

2. {\it The necessity of the compactness of $i: X_1\to X_2$ for conditions 
a) 
and b) to hold.}

  Assume now that  $i$ is compact. Let us prove that 
conditions a) and b) hold. In the proof of the necessity
of these conditions the assumption about the
compatibility of the norms of $X_2$ and $X_3$ is used essentially.
Without this assumption one cannot prove that conditions a) and b) hold.
This is proved in the {\bf Remark 1} after the end of the  proof of 
Theorem 1.
 
If  $i$ is compact, then assumption a) holds because  $||u||_2\geq 
||u||_3$. 
Suppose that assumption b) fails.
Then there is a sequence $u_n$ and a number $s_0>0$ such 
that $||u_n||_1=1$ and
$$
||u_n||_2\geq s_0+n||u_n||_3.
\eqno{(1)}$$
If  the embedding operator $i$ is compact and $||u_n||_1=1$, 
then one may assume that
the sequence $u_n$ converges in $X_2$. Its limit cannot be equal to zero, 
because, by (1), $||u_n||_2\geq s_0>0$. The  sequence $u_n$ converges in 
$X_3$  because $||u_n-u_m||_2\geq ||u_n-u_m||_3$, and {\it its 
limit in $X_3$ is not
zero, because the norms in $X_3$ and in $X_2$ are compatible.} Thus,
$\lim_{n\to \infty}||u_n||_3>0$.

Thus, (1) implies $||u_n||_3=O(\frac 1 n)\to 0$ as $n\to \infty$,
while $\lim_{n\to \infty}||u_n||_3>0$. This is a contradiction, which 
proves that b) holds.

Theorem 1 is proved. \hfill $\Box$

{\bf Remark 1.} In [1], p. 35, under the name Lions' lemma, the following 
claim is stated: 

{\bf Claim} ([1], p.35): {\it Let $X_1\subset X_2 \subset X_3$ be three 
Banach 
spaces.
Suppose the embedding $X_1\to X_2$ is compact. Then given any 
$\epsilon>0$, there is a $K(\epsilon)>0$, such that 
$||u||_2\leq \epsilon ||u||_1+K(\epsilon)||u||_3$ for all $u\in X_1$.}

This claim, is {\it not correct} because there is no assumption about
compatibility of the norms of $X_2$ and $X_3$.

For example, 
let $L^2(0,1)$ be the usual Lebesgue 
space of square integrable functions, $X_3=L^2(0,1),$
and $X_2$ be a Banach space of $L^2(0,1)$ functions with a finite
value at a fixed point $y\in [0,1]$ and with the 
norm  
$$||u||_2:=||u||_{L^2(0,1)}+|u(y)|=||u||_3+|u(y)|.$$ 
The space $X_2$ is complete
because $X_3$ is complete and the one-dimensional space, consisting of
numbers $u(y)$ with the usual norm $|u(y)|$, is complete. A   function 
$u_0(x)=0$ for $x\neq 0$ and $u_0(y)=1$ has the properties 
$$||u_0||_3=0,\quad ||u_0||_2=1.$$
Clearly,  $X_2\subset X_3$, and the norms in $X_2$ and $X_3$ are 
{\it comparable}, i.e., $||u||_3\leq ||u||_2$. 
However, these norms  
are {\it not compatible}: there is a convergent to zero 
sequence $\lim_{n\to \infty}u_n= 0$ in $X_3$ such that it does not 
converge 
to zero in  $X_2$, for example,
$\lim_{n\to \infty}||u_n||_2= 1$ in $X_2$. For instance, one may take 
$u_n(x)=u_0(x)$ for all $n=1,2,\dots $,
and an arbitrary fixed $y\in [0,1]$. Then $||u_n||_2=1$ and 
$||u_n||_3=0$, $\lim_{n\to \infty}||u_n||_2=1$ and  $\lim_{n\to 
\infty}||u_n||_3=0$. The sequence $u_n$ convereges to zero in $X_3$
and to a non-zero element $u_0$ in $X_2$.
In this case inequality (1) holds for any 
fixed $s_0\in (0,1)$ and any $n$, but the 
contradiction, which was used in the proof of the necessity in Theorem 1, 
can not be obtained  because $||u_n||_3=0$ for all $n$. 

Let us construct a  counterexample which shows that 
the Claim in [1], mentioned above, is not correct. Fix a $y\in [0,1]$.
Choose  the one-dimensional space of functions $\{u: 
u=\lambda u_0(x)\}$ as $X_1$, where 
$\lambda=const$, and define the norm in $X_1$ by the formula 
$||u||_1=|\lambda|$. Let $X_3=L^2(0,1)$. 
The space $X_1$ is a one-dimensional Banach space. Therefore bounded sets
in $X_1$ are precompact. Note that $|\lambda|=||\lambda u_0||_1=||\lambda 
u_0||_2\geq ||\lambda u_0||_3=0$ because $||u_0||_3=0$. Here the Banach 
space $X_2$ is defined as above with the norm  
$||u||_2:=||u||_{L^2(0,1)}+|u(y)|$, and the equalities
$||u_0||_2=1$ and $||u_0||_3=0$ are used. 

Consequently, 
$$X_1\subset X_2\subset X_3, \quad ||u||_1\geq ||u||_2\geq 
||u||_3,$$ 
and the embedding
$i: X_1\to X_2$ is compact because bounded sets in 
finite-dimensional 
spaces are precompact and $X_1$ is a one-dimensional space. Thus, all the 
assumptions 
of the  {\bf Claim} are satisfied.
However the inequality of the {\bf Claim}:
$$||u||_2\leq \epsilon ||u||_1+K(\epsilon)||u||_3 \quad \forall  u\in 
X_1$$
does not hold for any fixed $\epsilon \in (0,1)$. In our 
counterexample $u=\lambda u_0$, $||u_0||_3=0$,
and the above inequality takes the 
form: $|\lambda|\leq \epsilon |\lambda|$.
Clearly, this inequality does not hold 
for a fixed $\epsilon \in (0,1)$ unless $\lambda =0$.

\end{document}